\documentclass[twoside,reqno]{amsart}

\usepackage{amsmath,amssymb,amsbsy,amstext,amsxtra,latexsym}
\usepackage{graphicx}
\usepackage{subfig}
\usepackage{enumerate}
\usepackage[all]{xy}

\newtheorem{theorem}{Theorem}[section]

\newtheorem{proposition}[theorem]{Proposition}

\theoremstyle{definition}

\theoremstyle{remark}

\numberwithin{equation}{section}

\newcommand{\uc}[1]{\ensuremath \overset{#1}{\circ}}

\newcommand{\blup}[2]{\ensuremath #1 \sharp #2 \overline{\mathbb{CP}^2}}

\def\sheaf#1{\ensuremath \mathcal#1}

\begin{document}

\title[A surface of general type with $p_g=1$, $q=0$, and $K^2=8$]{A simply connected surface of general type \\
with $p_g=1$, $q=0$, and $K^2=8$}

\author{Heesang Park}

\address{Department of Mathematical Sciences, Seoul National University, 599 Gwanak-ro, Gwanak-gu, Seoul 151-747, Korea}

\email{hspark@math.snu.ac.kr}

\author{Jongil Park}

\address{Department of Mathematical Sciences, Seoul National University, 599 Gwanak-ro, Gwanak-gu, Seoul 151-747, Korea}

\email{jipark@snu.ac.kr}

\author{Dongsoo Shin}

\address{Department of Mathematics, Chungnam National University, Daejeon 305-764, Korea}

\email{dsshin@cnu.ac.kr}

\date{October 19, 2009}

\subjclass[2000]{Primary 14J29; Secondary 14J10, 14J17, 53D05}

\keywords{$\mathbb{Q}$-Gorenstein smoothing, rational blow-down surgery, surface of general type}

\begin{abstract}
We construct a new family of simply connected minimal complex surfaces with $p_g=1$, $q=0$, and $K^2=8$ using a $\mathbb{Q}$-Gorenstein smoothing theory.
\end{abstract}

\maketitle

\section{Introduction}

This paper is an addendum to the authors' work~\cite{PPS-p_g=1}, in which we constructed a family of minimal complex surfaces of general type with $p_g=1$, $q=0$, and $1 \le K^2 \le 2$ and simply connected surfaces with $p_g=1$, $q=0$, and $3 \le K^2 \le 6$ using a $\mathbb{Q}$-Gorenstein smoothing theory. We extend the results to the $K^2=8$ case in this paper. The main result of this paper is the following theorem.

\begin{theorem}
\label{thm-main}
There exists a simply connected minimal complex surface of general type with $p_g=1$, $q=0$, and $K^2=8$.
\end{theorem}

We briefly sketch the proof. At first we blow up a K3 surface $\overline{Y}$ in a suitable set of points so that we obtain a surface with some special disjoint linear chains of rational curves which can be contracted to singularities class $T$ on a singular surface $\overline{X}$ with  $H^2(\sheaf{T_{\overline{X}}}) \neq 0$. In order to prove the existence of a global $\mathbb{Q}$-Gorenstein smoothing of $\overline{X}$, we apply the cyclic covering trick developed in Y. Lee and J. Park~\cite{Lee-Park-Horikawa}. The cyclic covering trick says that, if a cyclic covering $\pi : V \to W$ of singular surfaces satisfies certain conditions and the base $W$ has a $\mathbb{Q}$-Gorenstein smoothing, then the cover $V$ has also a $\mathbb{Q}$-Gorenstein smoothing.

The main ingredient of this paper is that we construct an unramified double covering $\overline{\pi} : \overline{X} \to X$ to a singular surface $X$ constructed in a recent paper \cite{Park} of the first author. It is a main result of H. Park~\cite{Park} that the singular surface $X$ has a global $\mathbb{Q}$-Gorenstein smoothing and a general fiber $X_t$ of the smoothing of $X$ is a surface of general type with $p_g=0$, $K^2=4$, and $\pi_1 = \mathbb{Z}/2\mathbb{Z}$. We show that the double covering $\overline{X} \to X$ satisfies all the conditions of the cyclic covering trick; hence, there is a global $\mathbb{Q}$-Gorenstein smoothing of $\overline{X}$. Then it is not difficult to show that a general fiber $\overline{X}_t$ of the smoothing of $\overline{X}$ is the desired surface.

\section{Construction}

\begin{figure}[hbtb]
\centering
\includegraphics{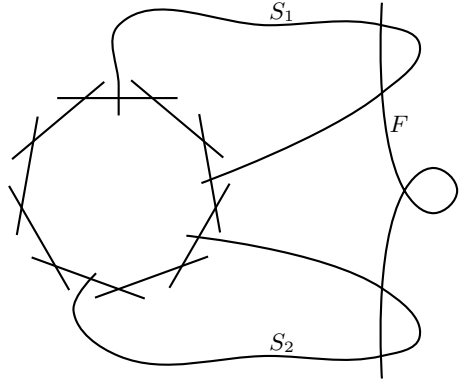}
\caption{An Enriques surface $Y$}
\label{figure:Y}
\end{figure}

\begin{figure}[hbtb]
\centering
\includegraphics{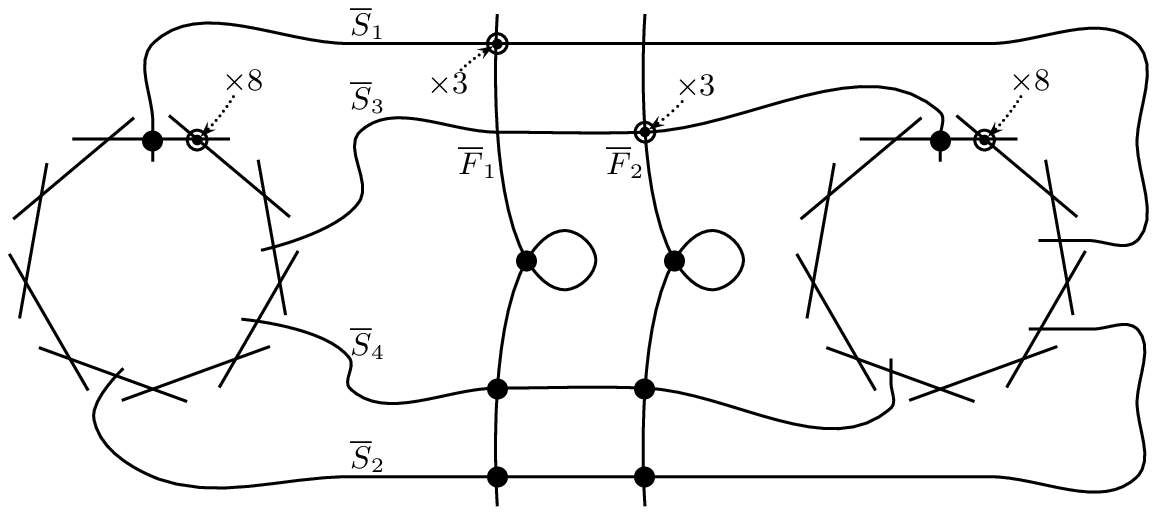}
\caption{A K3 surface $\overline{Y}$}
\label{figure:Ybar}
\end{figure}

According to Kondo~\cite{Kondo}, there is an Enriques surface $Y$ with an elliptic fibration over $\mathbb{P}^1$ which has a $I_9$-singular fiber, a nodal singular fiber $F$, and two bisections $S_1$ and $S_2$; Figure~\ref{figure:Y}. Again by Kondo~\cite{Kondo}, there is an unbranched double covering $\pi: \overline{Y} \to Y$ of $Y$ where $\overline{Y}$ is an elliptic K3 surface  which has two $I_9$-singular fiber, two nodal singular fiber $\overline{F}_1$ and $\overline{F}_2$, and four sections $\overline{S}_1, \dotsc, \overline{S}_4$ such that $\pi(\overline{F}_1) = \pi(\overline{F}_2) = F$, $\pi(\overline{S}_1) = \pi(\overline{S}_3) = S_1$, and $\pi(\overline{S}_2) = \pi(\overline{S}_4) = S_2$; Figure~\ref{figure:Ybar}.

\begin{figure}[hbtb]
\centering
\includegraphics{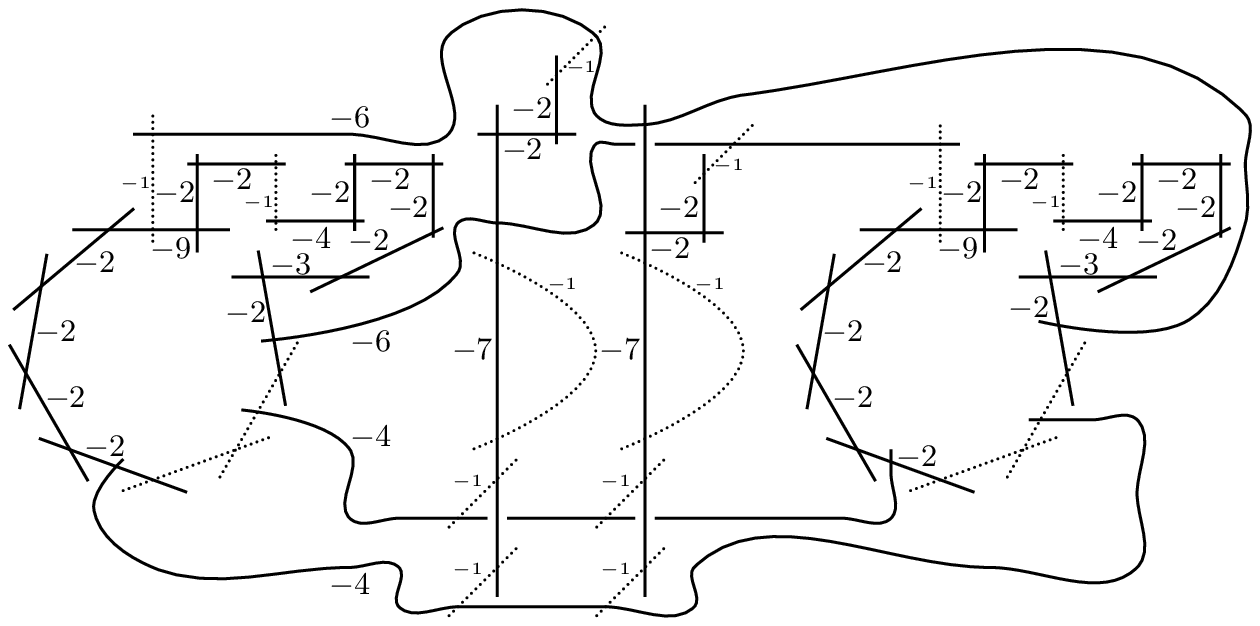}
\caption{A surface $\overline{Z} = \blup{\overline{Y}}{30}$}
\label{figure:Zbar}
\end{figure}

We blow up the K3 surface $\overline{Y}$ totally 30 times at the marked points $\bullet$ and $\bigodot$. We then get a surface $\overline{Z} = \blup{\overline{Y}}{30}$; Figure~\ref{figure:Zbar}. There exist four disjoint linear chains of $\mathbb{CP}^1$'s in $\overline{Z}$:
\begin{align*}
&C_{19,6}: \uc{-2}-\uc{-2}-\uc{-9}-\uc{-2}-\uc{-2}-\uc{-2}-\uc{-2}-\uc{-4} \\
&C_{19,6}: \uc{-2}-\uc{-2}-\uc{-9}-\uc{-2}-\uc{-2}-\uc{-2}-\uc{-2}-\uc{-4} \\
&C_{73,50}: \uc{-2}-\uc{-2}-\uc{-7}-\uc{-6}-\uc{-2}-\uc{-3}-\uc{-2}-\uc{-2}-\uc{-2}-\uc{-2}-\uc{-4}\\
&C_{73,50}: \uc{-2}-\uc{-2}-\uc{-7}-\uc{-6}-\uc{-2}-\uc{-3}-\uc{-2}-\uc{-2}-\uc{-2}-\uc{-2}-\uc{-4}
\end{align*}

We contract these four chains of $\mathbb{CP}^1$'s from the surface $\overline{Z}$ so that it produces a normal projective surface $\overline{X}$ with four singular points of class $T$. It is not difficult to show that $H^2(X, \sheaf{T_X}) \neq 0$.

\begin{theorem}\label{theorem:Q-Gorenstein}
The singular surface $\overline{X}$ has a global $\mathbb{Q}$-Gorenstein smoothing. A general fiber $\overline{X}_t$ of the smoothing of $\overline{X}$ is a simply connected minimal complex surface of general type with $p_g=1$, $q=0$, and $K^2=8$.
\end{theorem}

In order to prove Theorem~\ref{theorem:Q-Gorenstein}, we apply the following proposition.

\begin{proposition}[{Y. Lee and J. Park~\cite{Lee-Park-Horikawa}}]\label{proposition:Lee-Park}
Let $V$ be a normal projective surface with singularities of class $T$. Assume that a cyclic group $G$ acts on $X$ such that
\begin{enumerate}[1.]
\item $W = V/G$ is a normal projective surface with singularities of class $T$,

\item $p_g(W)=q(W)=0$,

\item $W$ has a $\mathbb{Q}$-Gorenstein smoothing,

\item the map $\sigma : V \to W$ induced by a cyclic covering is flat, and the branch locus $D$ (resp. the ramification locus) of the map $\sigma : V \to W$ is a nonsingular curve lying outside the singular locus of $W$ (resp. of $V$), and

\item $H^1(W, \sheaf{O_W}(D))=0$.
\end{enumerate}
Then there exists a $\mathbb{Q}$-Gorenstein smoothing of $V$ that is compatible with a $\mathbb{Q}$-Gorenstein smoothing of $W$. Furthermore the cyclic covering extends to the $\mathbb{Q}$-Gorenstein smoothing.
\end{proposition}

\begin{figure}[hbtb]
\centering
\includegraphics{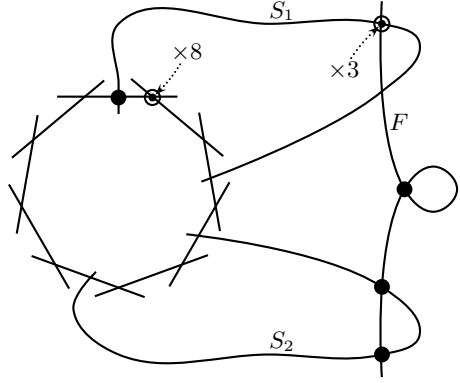}
\caption{An Enriques surface $Y$ with marked points}
\label{figure:Y-dots}
\end{figure}

\begin{figure}[hbtb]
\centering
\includegraphics{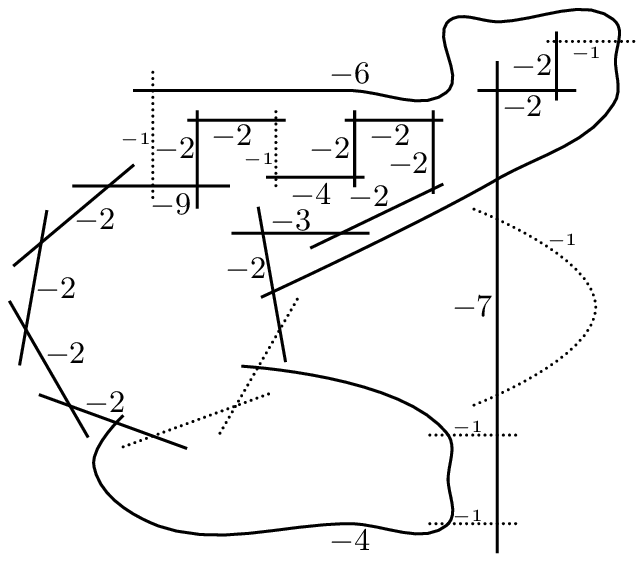}
\caption{A surface $Z = \blup{Y}{15}$}
\label{figure:Z}
\end{figure}

We now construct an unramified double covering from the singular surface $\overline{X}$ to another singular surface. We begin with the Enriques surface $Y$ in Figure~\ref{figure:Y}. We blow up totally 15 times at the marked points $\bullet$ and $\bigodot$; Figure~\ref{figure:Y-dots}. We then get a surface $Z = \blup{Y}{15}$; Figure~\ref{figure:Z}. There exist two disjoint linear chains of $\mathbb{CP}^1$'s in $Z$:
\begin{align*}
&C_{19,6}: \uc{-2}-\uc{-2}-\uc{-9}-\uc{-2}-\uc{-2}-\uc{-2}-\uc{-2}-\uc{-4} \\
&C_{73,50}: \uc{-2}-\uc{-2}-\uc{-7}-\uc{-6}-\uc{-2}-\uc{-3}-\uc{-2}-\uc{-2}-\uc{-2}-\uc{-2}-\uc{-4}
\end{align*}

We contract the two chains of $\mathbb{CP}^1$'s from the surface $Z$ so that it produces a normal projective surface $X$ with two singular points class $T$. It is clear that there is an unbranched double covering $\overline{\pi} : \overline{X} \to X$. The singular surface $X$ satisfies the third condition of Proposition~\ref{proposition:Lee-Park}.

\begin{proposition}[{H. Park~\cite{Park}}]\label{proposition:Park}
The singular surface $X$ has a global $\mathbb{Q}$-Gorenstein smoothing. A general fiber $X_t$ of the smoothing of $X$ is a minimal complex surface of general type with $p_g=0$, $K^2=4$, and $\pi_1(X_t) = \mathbb{Z}/2\mathbb{Z}$.
\end{proposition}

\begin{proof}[Proof of Theorem~\ref{theorem:Q-Gorenstein}]
It is easy to show that the covering $\overline{\pi} : \overline{X} \to X$ satisfies all conditions of Proposition~\ref{proposition:Lee-Park}. Therefore the singular surface $\overline{X}$ has a global $\mathbb{Q}$-Gorenstein smoothing. Let $\overline{X}_t$ be a general fiber of the smoothing of $\overline{X}$. Since $p_g(\overline{X})=1$, $q(\overline{X})=0$, and $K_{\overline{X}}^2=8$, by applying general results of complex surface theory and $\mathbb{Q}$-Gorenstein smoothing theory, one may conclude that a general fiber $\overline{X}_t$ is a complex surface of general type with $p_g=1$, $q=0$, and $K^2=8$. Furthermore, it is not difficult to show that a general fiber $X_t$ is minimal by using a similar technique in~\cite{Lee-Park-K^2=2, PPS-K3, PPS-K4}.

\textit{Claim.} A general fiber $\overline{X}_t$ is simply connected: By Proposition~\ref{proposition:Lee-Park}, there is an induced unbranched double covering $\overline{X}_t \to X_t$; hence, a general fiber $\overline{X}_t$ is simply connected because $\pi(X_t) = \mathbb{Z}/2\mathbb{Z}$ by Proposition~\ref{proposition:Park}.
\end{proof}

\end{document}